\input amstex


\def\b1{\text{\bf 1}}

\def\gr{\text{gr}}

\def\Hom{\text{Hom}}

\def\#{\,\check{}}

\def\id{\text{id}}

\def\Rep{\text{Rep}}

\def\Gal{\text{Gal}}

\def\st{\text{st}}


\def\hra{\hookrightarrow}
\def\iso{\buildrel\sim\over\rightarrow}


\documentstyle{amsppt}
\NoBlackBoxes

\topmatter
\title       On a theorem of Kisin       \endtitle
\author Alexander Beilinson and Floric Tavares Ribeiro \endauthor
\leftheadtext{A.~Beilinson and F.~Tavares Ribeiro}

\address A.B.: Department of Mathematics, University of Chicago, Chicago, IL 60637, USA,  \newline F.T.R.: 
\'Ecole Normale Sup\'eriore de Lyon, Unit\'e de Math\'ematiques Pures et Appliqu\'ees, 46, All\'ee d'Italie, 69364 LYON Cedex 07, France
\endaddress

\email sasha\@math.uchicago.edu, floric.tavares-ribeiro@unicaen.fr \endemail

\thanks The first named author was supported in part by NSF
Grant DMS-0401164.\endthanks




\endtopmatter

\document


Let $K$ be a $p$-adic field, i.e., a complete discretely valued field of characteristic 0 with perfect residue field of characteristic $p>0$, and $\bar{K}$ be an algebraic closure of $K$. We fix a uniformiser $\pi\in K$. Let $\Xi =\Xi_\pi$ be the corresponding Kummer $\Bbb Z_p (1)$-torsor; its  elements are sequences
$\xi = (\xi_n )_{n\ge 0}$ of elements in $\bar{K}$ such that $\xi_{n+1}^p =\xi_n$, $\xi_0 =\pi$. Pick one $\xi$, and set $K_\xi =\cup K(\xi_n )$. Consider the Galois groups
 $G :=\Gal (\bar{K}/K)$, $G_{\xi }:=\Gal (\bar{K}/K_\xi )$; let $\Rep (G)$, $\Rep (G_{\xi})$ be the categories of their finite-dimensional $\Bbb Q_p$-representations.
 
 The next result was conjectured by Breuil \cite {B} and
 proved by Kisin  \cite{K} 0.2; the proof in loc.~cit.~is based on  theory of Kisin modules. This note provides an alternative argument that uses only basic properties of Fontaine's rings; its key ingredient  (namely, (i) of the lemma below) is the same as in Grothendieck's proof of the monodromy theorem.
 
 \proclaim{Theorem} 
 The restriction functor $\Rep (G)\to \Rep (G_{\xi})  $ is fully faithful on
 the subcategory of crystalline representations. 
 \endproclaim

 \demo{Proof}  The Galois group $G$ acts on $\Xi$, and $G_\xi$ is the stabilizer of $\xi$. The action is transitive, i.e., $G/G_\xi \iso \Xi$, since polynomials $t^{p^n}-\pi$ are irreducible.
 
   Let $R$ be the ring of  continuous $\Bbb Q_p$-valued functions on  $\Xi$. Let  $ R_{\st}\subset R_\phi$ be the subrings of polynomial, resp.~locally polynomial,  functions (this makes sense since $\Xi$ is $\Bbb Z_p (1)$-torsor). Since $G$ acts on $\Xi$ by affine transformations, its action on $R$ preserves the subrings.
   
    \proclaim{Lemma} 
(i) $R_\phi$ is the union of all finite-dimensional $G$-submodules of $R$. \newline (ii)  $R_{\st}$ is the union of all semi-stable $G$-submodules of $R_\phi$. \newline (iii) $\Bbb Q_p$ is the only nontrivial crystalline $G$-submodule of $R_{\st}$.
 \endproclaim
   
Assuming the lemma, let us prove the theorem.  For $V\in \Rep (G_\xi )$ we denote by $I(V)$ the induced $G$-module, that is the space of all continuous maps $f: G \to V$ such that $f(hg)=hf(g)$ for $h\in G_{\xi}$, the action of $G$ is $g(f)(g')=f(g'g)$. It is a $G$-equivariant $R$-module, the $R$-action is $(rf)(g)=r(g^{-1}\xi)f(g)$.
For $U\in \Rep (G)$ we have the Frobenius reciprocity  $ \Hom_{G_\xi}(U,V) \iso \Hom_G (U, I(V))$ that identifies $\alpha : U\to V$ with $\tilde{\alpha}: U\to I(V)$, $\tilde{\alpha}(u)(g)=\alpha (gu)$,  $\alpha (u)=\tilde{\alpha}(u)(1)$. For   $V\in\Rep (G)$ the image of $\id_V \in  \Hom_{G_\xi}(V,V)$ is a $G$-morphism $V\to I(V)$ that yields an identification of $G$-equivariant $R$-modules $V\otimes R  \iso I(V)$.

So for $V_1 ,V_2\in\Rep (G)$ one has identifications $\Hom_{G_\xi}(V_1 ,V_2 )= \Hom_G (V_1,I(V_2))$ $=\Hom_G (V_1 ,V_2 \otimes R) = \Hom_G (V_1\otimes V_2^* ,R)=  \Hom_G (V_1 \otimes V_2^* ,R_{\phi} )$,  the last equality comes from (i). 
 If both  $V_i$ are crystalline, then this equals $ \Hom_G (V_1\otimes V_2^* ,\Bbb Q_p )=\Hom_G (V_1 ,V_2 )$ by (ii), (iii). Thus $\Hom_{G_\xi}(V_1 ,V_2 )=\Hom_G (V_1 ,V_2 )$, q.e.d. \qed
   
  \demo{Proof of Lemma} Let $P$ be the group of all  affine automorphisms of $\Bbb Z_p (1)$-torsor $\Xi$; it is an extension of $\Bbb Z_p^\times $ by $\Bbb Z_p (1)$,  the choice of $\xi$ gives a splitting. Let $\eta: G \to P$ be the action of $G$ on $\Xi$; its composition with 
 $P\twoheadrightarrow \Bbb Z_p^\times $ is   cyclotomic character $\chi$. 
 
 Consider the filtration  $R_{\st\, n}$ on $R_{\st}$ by the degree of the polynomial. Then $G$ acts on $\gr_n R_{\st}$ by $\chi^{-n}$, i.e., $\gr_n R_{\st}$ is isomorphic to $\Bbb Q_p (-n)$. 
 
 There is a  canonical morphism $\varepsilon : R_{\st}\to$ B$_{\st}$ of $\Bbb Q_p$-algebras defined as follows. For $\xi \in \Xi$ let $l_\xi : \Xi \to \Bbb Z_p (1)$ be the identification of torsors such that $l_\xi (\xi )=0$. If $\tau$ is a generator of $\Bbb Z_p (1)$, then $\tau^{-1}l_\xi \in R_{\st}$ is a linear polynomial function, i.e., a free generator of $R_{\st}$. We define $\varepsilon$ by formula $\varepsilon ( \tau^{-1}l_\xi)= - \tau^{-1}\lambda(\xi)$. Here in the r.h.s.~we view $\tau$ as an invertible element of B$_{\text{crys}}$ via the embedding $\Bbb Z_p (1)\hra $ B$_{\text{crys}}$ from \cite{F1} 2.3.4, and $\lambda (\xi )\in $ B$_{\st}$ is as in \cite{F1} 3.1.4. It follows from the definitions in \cite{F1} 3.1 that  $\varepsilon$ does not depend on the auxiliary choice of $\xi$. It evidently commutes with the Galois action.
 Since  $\log (\xi ) $ is a free generator of B$_{\st}$ over B$_{\text{crys}}$, we see that $\varepsilon$ is injective and   $R_{\st\, n}$ for $n\ge 1$ are {\it non-crystalline} semi-stable $G$-modules. 
 
 Choose  $v$ and $\log$ from \cite{F2} 5.1.2 as $v(\pi )=1$, $\log (\pi )=0$. As in \cite{F2} 5.2, this yields the fully faithful  tensor functor $D_{\st }:\Rep(G)_{\st}\to \text{MF}_K(\varphi ,N)$. 
 
 Consider the polynomial algebra $K_0 [t]$. We equip it with Frobenius semi-linear automorphism $\varphi$, $\varphi (t):=pt$, the $K_0$-derivation $N:=\partial_t$, and the Hodge filtration  $F^i$:= the $K$-span of $t^{\ge i}$. The subspaces of polynomials of degree $\le n$ are filtered $(\varphi ,N)$-modules,  so $K_0 [t]$ is a ring ind-object of MF$_K (\varphi ,N)$.
 
 There is a canonical isomorphism $K_0 [t] \iso D_{\st }(R_{\st})$ which identifies $t$ with $(\tau^{-1}l_\xi )\otimes \tau + 1\otimes\lambda (\xi ) \in ( R_{\st}\otimes \text{B}_{\st})^G = D_{\st }(R_{\st})$. Thus each $D_{\st}(R_{\st\, n})$ is a single Jordan block for the action of $N$, so every finite-dimensional $G$-submodule of $R_{\st}$ equals one of $R_{\st\, n}$'s, which implies (iii).

Notice that $R_\phi = R_0 \otimes R_{\st}$, where $R_0$ is the subring of locally constant functions. Since $G$ acts transitively on $\Xi$, one has $R_0^{G}=\Bbb Q_p$ and all finite-dimensional  $G$-modules that occur in $R_0$ are generated by $G_{\xi}$-fixed vectors. These representations are Artinian, hence semisimple, so we have the decomposition $R_0 =\Bbb Q_p \oplus R'_0$. Since the map $G_{\xi} \to\Gal (K^{\text{un}}/K)$, where $K^{\text{un}}\subset \bar{K}$ is the maximal unramified extension of $K$, is surjective (for $K^{\text{un}}\cap K_\xi =K$), every $G$-module in $R'_0$ is ramified. Thus every irreducible subquotient of $R'_0 \otimes R_{\st}$ is {\it not} semi-stable, and we get (ii).

It remains to prove (i). We first show that $\eta (G)$ is an open subgroup of $P$. Since $\chi (G)$ is  open in $\Bbb Z_p^\times$, it suffices to check that $\eta (G)\cap \Bbb Z_p (1)$ is open in $\Bbb Z_p (1)$.
Since every closed nontrivial subgroup of $\Bbb Z_p (1)$ is open, we need to check that $\eta (G)\cap \Bbb Z_p (1)\neq \{ 0\}$. If not, then $\eta (G)\iso \chi (G)$ is commutative, so $G$ acts on $R$ through an abelian quotient. This implies, since $\gr_n R_{\st}\simeq \Bbb Q_p (-n)$ are pairwise non-isomorphic $G$-modules,  that filtration $R_{\st\, n}$ splits, which is not true, q.e.d.

 Let $\tau$ be a
generator of $\Bbb Z_p (1)\subset P$; then $R_\phi$ is the union of all finite-dimensional $   \Bbb Z_p (1)$-submodules of $R$ on which all eigenvalues of $\tau$ 
are roots of 1. Since $\eta (G)$ has finite index in $P$, it suffices to show that every finite-dimensional $P$-submodule $V$ of $R$ has this property. This follows since
  for $g\in P$ one has $g\tau g^{-1}= \tau^m$, where $m$ is the image of $g$ in $\Bbb Z_p^\times$, and there are only finitely many eigenvalues of $\tau$ on $V$. \qed
  \enddemo 
   
    \enddemo                                                             
                                                               
The above proof  was first found by F.T.R.~and later and independently by A.B.
We are grateful to Bhargav Bhatt, Matt Emerton, and Mark Kisin (A.B.), and to
 Laurent Berger and Fran\c{c}ois Brunault (F.T.R.) for very helpful discussions.

\end